\documentclass[12pt]{article}

\usepackage{amsmath, epsfig, cite}
\usepackage{amssymb}
\usepackage{amsfonts}
\usepackage{latexsym}
\usepackage{float}
\usepackage{color}

\newtheorem{thm}{Theorem}[section]

\newtheorem{conj}[thm]{Conjecture}
\newtheorem{lem}[thm]{Lemma}

\numberwithin{equation}{section}

\newcommand{\qed}{{\hfill$\square$}\medskip}

\begin{document}

\begin{center}
{\large\bf Some supercongruences arising from symbolic summation}
\end{center}

\vskip 2mm \centerline{Ji-Cai Liu}
\begin{center}
{\footnotesize Department of Mathematics, Wenzhou University, Wenzhou 325035, PR China\\
{\tt jcliu2016@gmail.com} }
\end{center}


\vskip 0.7cm \noindent{\bf Abstract.}
Based on some combinatorial identities arising from symbolic summation, we extend two supercongruences on partial sums of hypergeometric series, which were originally conjectured by Guo and Schlosser and recently confirmed by Jana and Kalita.


\vskip 2mm
\noindent{\it MR Subject Classifications}: 11A07, 05A19, 33F10, 33C20

\section{Introduction}
In 1997, Van Hamme \cite[(A.2)--(M.2) ]{van-b-1997} conjectured 13 supercongruences which relate partial sums of certain hypergeometric series to the values of $p$-adic gamma functions. Van Hamme's (D.2) supercongruence states that for primes $p\equiv 1\pmod{6}$,
\begin{align}
\sum_{k=0}^{p-1}(6k+1)\frac{\left(\frac{1}{3}\right)_k^6}{(1)_k^6}\equiv
 -p\Gamma_p\left(\frac{1}{3}\right)^9 \pmod{p^4}.\label{new-2}
\end{align}
Here the Pochhammer's symbol is given by
\begin{align*}
(x)_0=1\quad\text{and}\quad (x)_k=x(x+1)\cdots(x+k-1),
\end{align*}
and the $p$-adic gamma function is defined as
\begin{align*}
\Gamma_p(x)=\lim_{m\to x}(-1)^m\prod_{\substack{0< k < m\\
(k,p)=1}}k,
\end{align*}
where the limit is for $m$ tending to $x$ $p$-adically in $\mathbb{Z}_{\ge 0}$.

By using the Dougall's formula, Long and Ramakrishna \cite[Theorem 2]{lr-am-2016} generalized
Van Hamme's supercongruence \eqref{new-2} as follows:
\begin{align}
\sum_{k=0}^{p-1}(6k+1)\frac{\left(\frac{1}{3}\right)_k^6}{(1)_k^6}
\equiv \begin{cases}
 -p\Gamma_p\left(\frac{1}{3}\right)^9 \pmod{p^6}\quad &\text{if $p\equiv 1\pmod{6}$,}\\[10pt]
-\frac{10}{27}p^4\Gamma_p\left(\frac{1}{3}\right)^9 \pmod{p^6}\quad &\text{if $p\equiv 5\pmod{6}$.}
\end{cases}\label{new-3}
\end{align}

Recently, Guo and Schlosser \cite[Theorem 2.3]{gs-ca-2019} established a partial $q$-analogue of
\eqref{new-3} by using transformation formulas for basic hypergeometric series. One can refer to
\cite{guo-jmaa-2018,gl-jdea-2018,gs-rm-2019,gz-am-2019} for more $q$-analogues of congruences ($q$-congruence) for indefinite sums of binomial coefficients as well as hypergeometric series.

Guo and Schlosser \cite{gs-ca-2019} also proposed two conjectural supercongruences related to \eqref{new-3}.
\begin{conj}(See \cite[Conjecture 12.6]{gs-ca-2019}.)
For primes $p\equiv 1\pmod{3}$, we have
\begin{align}
\sum_{k=0}^{(p-1)/3}(6k+1)\frac{\left(\frac{1}{3}\right)_k^4(1)_{2k}}{(1)_k^4\left(\frac{2}{3}\right)_{2k}}\equiv p\pmod{p^3}.\label{a-1}
\end{align}
\end{conj}

\begin{conj}(See \cite[Conjecture 12.7]{gs-ca-2019}.)
For primes $p\equiv 2\pmod{3}$, we have
\begin{align}
\sum_{k=0}^{(p+1)/3}(6k-1)\frac{\left(-\frac{1}{3}\right)_k^4(1)_{2k}}{(1)_k^4\left(-\frac{2}{3}\right)_{2k}}\equiv p\pmod{p^3}.\label{a-2}
\end{align}
\end{conj}

Supercongruences \eqref{a-1} and \eqref{a-2} were first confirmed by Jana and Kalita \cite{jk-rj-2019} by using hypergeometric series identities and $p$-adic gamma functions.

The motivation of this paper is to extend \eqref{a-1} and \eqref{a-2} to strong versions modulo
$p^4$ (Theorems \ref{t-1} and \ref{t-2}). In this paper, we will provide a different approach to \eqref{a-1} and \eqref{a-2}, which is based on combinatorial identities arising from symbolic summation.

\begin{thm}\label{t-1}
For odd primes $p\equiv 1\pmod{3}$, we have
\begin{align}
\sum_{k=0}^{(p-1)/3}(6k+1)\frac{\left(\frac{1}{3}\right)_k^4(1)_{2k}}{(1)_k^4\left(\frac{2}{3}\right)_{2k}}\equiv p+\frac{p^3}{9}B_{p-2}\left(\frac{1}{3}\right)\pmod{p^4},\label{a-3}
\end{align}
where the Bernoulli polynomials are given by
\begin{align*}
\frac{te^{xt}}{e^t-1}=\sum_{k=0}^{\infty}B_k(x)\frac{t^k}{k!}.
\end{align*}
\end{thm}

\begin{thm}\label{t-2}
For odd primes $p\equiv 2\pmod{3}$, we have
\begin{align}
\sum_{k=0}^{(p+1)/3}(6k-1)\frac{\left(-\frac{1}{3}\right)_k^4(1)_{2k}}{(1)_k^4\left(-\frac{2}{3}\right)_{2k}}\equiv p-p^3\left(\frac{1}{9}B_{p-2}\left(\frac{1}{3}\right)-2\right)\pmod{p^4}.\label{a-4}
\end{align}
\end{thm}
{\it Remark.} For some open supercongruence conjectures involving Bernoulli polynomials, one can refer to \cite[Conjecture 6.11]{sunzw-ffa-2017}.

The rest of this paper is organized as follows. Section 2 is devoted to establishing some congruences on harmonic sums and some combinatorial identities. We prove Theorems \ref{t-1} and \ref{t-2} in Sections 3 and 4, respectively.

\section{Preliminary results}
Let
\begin{align*}
&H_n^{(r)}=\sum_{k=1}^n\frac{1}{k^r},\\
&S_n^{(r)}=\sum_{k=1}^n \frac{1}{(3k-1)^r},\\
&T_n^{(r)}=\sum_{k=1}^n \frac{1}{(3k-2)^r},
\end{align*}
with the conventions that $H_n=H_n^{(1)}$, $S_n=S_n^{(1)}$ and $T_n=T_n^{(1)}$.
In order to prove Theorems \ref{t-1} and \ref{t-2}, we need some congruences on harmonic sums.
\begin{lem}
Let $p\ge 5$ be a prime. For $1\le k \le p-1$, we have
\begin{align}
&H_{p-1-k}-H_k\equiv pH_k^{(2)}\pmod{p^2},\label{b-1}\\
&H_{p-1-k}^{(2)}+H_k^{(2)}\equiv 0\pmod{p}.\label{b-2}
\end{align}
\end{lem}
{\it Proof.}
Note that for $1\le j\le p-1$,
\begin{align*}
\frac{1}{p-j}\equiv -\frac{1}{j}-\frac{p}{j^2}\pmod{p^2}.
\end{align*}
By the Wolstenholme's theorem \cite[page 114]{hw-b-2008}, we have
\begin{align*}
H_{p-1-k}=H_{p-1}-\sum_{j=1}^k\frac{1}{p-j}\equiv H_k+pH_k^{(2)}\pmod{p^2},
\end{align*}
which is equivalent to \eqref{b-1}.

Since for $1\le j\le p-1$,
\begin{align*}
\frac{1}{(p-j)^2}\equiv \frac{1}{j^2}\pmod{p},
\end{align*}
we have
\begin{align*}
H_{p-1-k}^{(2)}+H_k^{(2)}\equiv H_{p-1}^{(2)}\equiv 0\pmod{p},
\end{align*}
where we have utilized the Wolstenholme's theorem \cite[page 114]{hw-b-2008}.
\qed

\begin{lem}
Let $p\equiv 1\pmod{3}$ be a prime and $n=(p-1)/3$. Then
\begin{align}
&S_n\equiv 0\pmod{p},\label{b-3}\\
&S_n^{(2)}\equiv -\frac{1}{9}B_{p-2}\left(\frac{1}{3}\right)\pmod{p}.\label{b-4}
\end{align}
\end{lem}
{\it Proof.}
Since $n\equiv -1/3\pmod{p}$, we have
\begin{align*}
S_n=\sum_{k=1}^n\frac{1}{3k-1}=\frac{1}{3}\sum_{k=1}^n\frac{1}{k-1/3}\equiv \frac{1}{3}\sum_{k=1}^n\frac{1}{k+n}=\frac{1}{3}\left(H_{2n}-H_n\right)\pmod{p}.
\end{align*}
Then the proof of \eqref{b-3} follows from the above and \eqref{b-1}.

Note that
\begin{align*}
S_n^{(2)}=\sum_{k=1}^n\frac{1}{(3k-1)^2}=\frac{1}{9}\sum_{k=1}^n\frac{1}{(k-1/3)^2}
\equiv \frac{1}{9}\sum_{k=1}^n\frac{1}{(k+n)^2}
=\frac{1}{9}\left(H_{2n}^{(2)}-H_n^{(2)}\right)\pmod{p}.
\end{align*}
Combining \eqref{b-2}, the above congruence and the following congruence \cite[(9)]{lehmer-am-1938}:
\begin{align}
H_{\lfloor p/3\rfloor}^{(2)}\equiv \frac{1}{2}\left(\frac{-3}{p}\right)B_{p-2}\left(\frac{1}{3}\right)\pmod{p},\label{new-1}
\end{align}
where $\left(\frac{\cdot}{p}\right)$ denotes the Legendre symbol,
we reach the desired congruence \eqref{b-4}.
\qed

\begin{lem}
Let $p\equiv 2\pmod{3}$ be an odd prime and $n=(p+1)/3$. Then
\begin{align}
&T_n\equiv 0\pmod{p},\label{b-5}\\
&T_n^{(2)}\equiv \frac{1}{9}B_{p-2}\left(\frac{1}{3}\right)\pmod{p}.\label{b-6}
\end{align}
\end{lem}
{\it Proof.} The proofs of \eqref{b-5} and \eqref{b-6} are similar to those of \eqref{b-3}
and \eqref{b-4}, and we omit the details here.
\qed

We also need some combinatorial identities, which are discovered and proved by symbolic summation
package {\tt Sigma} developed by Schneider \cite{schneider-slc-2007}. One can refer to \cite{liu-jsc-2019} for the same approach to finding and proving identities of this type.

\begin{lem}\label{l-1}
For any positive integer $n$, we have
\begin{align}
&\sum_{k=0}^n(6k+1)\frac{\left(\frac{1}{3}\right)_k^2 (1)_{2k}(-n)_k\left(n+\frac{2}{3}\right)_k}
{(1)_k^2\left(\frac{2}{3}\right)_{2k} \left(n+\frac{4}{3}\right)_{k} \left(-n+\frac{2}{3}\right)_{k}}
=3n+1,\label{b-7}\\
&\sum_{k=0}^n(6k-1)\frac{\left(-\frac{1}{3}\right)_k^2 (1)_{2k}(-n)_k\left(n-\frac{2}{3}\right)_k}
{(1)_k^2\left(-\frac{2}{3}\right)_{2k} \left(-n+\frac{4}{3}\right)_{k} \left(n+\frac{2}{3}\right)_{k}}
=3n-1.\label{b-8}
\end{align}
\end{lem}

\begin{lem}\label{l-2}
For any positive integer $n$, we have
\begin{align}
&\sum_{k=0}^n(6k+1)\frac{\left(\frac{1}{3}\right)_k^2 (1)_{2k}(-n)_k\left(n+\frac{2}{3}\right)_k}
{(1)_k^2\left(\frac{2}{3}\right)_{2k} \left(n+\frac{4}{3}\right)_{k} \left(-n+\frac{2}{3}\right)_{k}}
\sum_{j=1}^k\left(\frac{1}{(3j)^2}-\frac{1}{(3j-2)^2}\right)\notag\\
&=-\frac{3n+1}{3}\left(-2H_n+3S_n^{(2)}+6S_n+2H_nS_n-3S_n^2\right)\notag\\
&+\frac{3n+1}{3}\left(2\sum_{k=1}^n\frac{H_k}{3k-1}+2\sum_{k=1}^n
\frac{T_k}{k}-6\sum_{k=1}^n\frac{T_k}{3k-1}\right),\label{b-9}
\end{align}
and
\begin{align}
&\sum_{k=0}^n(6k-1)\frac{\left(-\frac{1}{3}\right)_k^2 (1)_{2k}(-n)_k\left(n-\frac{2}{3}\right)_k}
{(1)_k^2\left(-\frac{2}{3}\right)_{2k} \left(-n+\frac{4}{3}\right)_{k} \left(n+\frac{2}{3}\right)_{k}}
\sum_{j=1}^k\left(\frac{1}{(3j)^2}-\frac{1}{(3j-4)^2}\right)\notag\\
&=\frac{3n-1}{3}\left(-\frac{(3n-1)(9n^2-8)}{n(3n-2)^2}+H_n-3T_n^{(2)}+3T_n^2-2H_nT_n+\frac{3(3n-4)T_n}{3n-2}\right)\notag\\
&+\frac{3n-1}{3}\left(-\frac{2(3n-1)S_{n-1}}{n}+2\sum_{k=1}^n\frac{H_k}{3k-2}+2\sum_{k=1}^{n-1}\frac{S_k}{k}-6\sum_{k=1}^{n-1}\frac{S_k}{3k-2}\right).
\label{b-10}
\end{align}
\end{lem}

\section{Proof of Theorem \ref{t-1}}
Letting $n=(p-1)/3$ in \eqref{b-7} gives
\begin{align}
\sum_{k=0}^{(p-1)/3}(6k+1)\frac{\left(\frac{1}{3}\right)_k^2 (1)_{2k}\left(\frac{1}{3}-\frac{p}{3}\right)_k\left(\frac{1}{3}+\frac{p}{3}\right)_k}
{(1)_k^2\left(\frac{2}{3}\right)_{2k} \left(1-\frac{p}{3}\right)_{k} \left(1+\frac{p}{3}\right)_{k}}=p.\label{c-1}
\end{align}
Note that
\begin{align*}
\frac{\left(\frac{1}{3}-\frac{p}{3}\right)_k\left(\frac{1}{3}+\frac{p}{3}\right)_k}
{\left(1-\frac{p}{3}\right)_{k} \left(1+\frac{p}{3}\right)_{k}}
=\prod_{j=1}^k\frac{(3j-2)^2-p^2}{(3j)^2-p^2}.
\end{align*}
From the following two Taylor expansions:
\begin{align*}
\frac{(3j-2)^2-x^2}{(3j)^2-x^2}=\left(\frac{3j-2}{3j}\right)^2-\frac{4(3j-1)}{(3j)^4}x^2+\mathcal{O}(x^4),
\end{align*}
and
\begin{align*}
\prod_{j=1}^k(a_j+b_jx^2)=\prod_{j=1}^ka_j\cdot\left(1+x^2\sum_{j=1}^k\frac{b_j}{a_j}\right)
+\mathcal{O}(x^4),
\end{align*}
we deduce that
\begin{align}
\frac{\left(\frac{1}{3}-\frac{p}{3}\right)_k\left(\frac{1}{3}+\frac{p}{3}\right)_k}
{\left(1-\frac{p}{3}\right)_{k} \left(1+\frac{p}{3}\right)_{k}}&\equiv \prod_{j=1}^k\left(\left(\frac{3j-2}{3j}\right)^2-\frac{4(3j-1)}{(3j)^4}p^2\right)\notag\\
&\equiv \frac{\left(\frac{1}{3}\right)_k^2}{(1)_k^2}
\left(1+p^2\sum_{j=1}^k\left(\frac{1}{(3j)^2}-\frac{1}{(3j-2)^2}\right)\right)\pmod{p^4}.
\label{c-2}
\end{align}
Substituting \eqref{c-2} into \eqref{c-1} gives
\begin{align}
&\sum_{k=0}^{(p-1)/3}(6k+1)\frac{\left(\frac{1}{3}\right)_k^4 (1)_{2k}}
{(1)_k^4\left(\frac{2}{3}\right)_{2k}}\notag\\
&\equiv p-p^2\sum_{k=0}^{(p-1)/3}(6k+1)\frac{\left(\frac{1}{3}\right)_k^4 (1)_{2k}}
{(1)_k^4\left(\frac{2}{3}\right)_{2k}}\sum_{j=1}^k\left(\frac{1}{(3j)^2}-\frac{1}{(3j-2)^2}\right)
\pmod{p^4}.\label{c-3}
\end{align}

Letting $n=(p-1)/3$ in \eqref{b-9} and noting \eqref{c-2}, we obtain
\begin{align*}
&\sum_{k=0}^{(p-1)/3}(6k+1)\frac{\left(\frac{1}{3}\right)_k^4 (1)_{2k}}
{(1)_k^4\left(\frac{2}{3}\right)_{2k}}\sum_{j=1}^k\left(\frac{1}{(3j)^2}-\frac{1}{(3j-2)^2}\right)\\
&\equiv -\frac{p}{3}\left(-2H_n+3S_n^{(2)}+6S_n+2H_nS_n-3S_n^2\right)\notag\\
&+\frac{p}{3}\left(2\sum_{k=1}^n\frac{H_k}{3k-1}+2\sum_{k=1}^n
\frac{T_k}{k}-6\sum_{k=1}^n\frac{T_k}{3k-1}\right)\pmod{p^2}.
\end{align*}
Combining \eqref{b-3} and the above, we obtain
\begin{align}
&\sum_{k=0}^{(p-1)/3}(6k+1)\frac{\left(\frac{1}{3}\right)_k^4 (1)_{2k}}
{(1)_k^4\left(\frac{2}{3}\right)_{2k}}\sum_{j=1}^k\left(\frac{1}{(3j)^2}-\frac{1}{(3j-2)^2}\right)\notag\\
&\equiv\frac{p}{3}\left(2\sum_{k=1}^n\frac{H_k}{3k-1}+2\sum_{k=1}^n
\frac{T_k}{k}-6\sum_{k=1}^n\frac{T_k}{3k-1}+2H_n-3S_n^{(2)}\right)\pmod{p^2}.\label{c-4}
\end{align}
By \eqref{b-1}, we have
\begin{align*}
&\sum_{k=1}^n\frac{H_k}{3k-1}\equiv\frac{1}{3}\sum_{k=1}^n\frac{1}{n+k}
\sum_{j=1}^k\frac{1}{j}=\frac{1}{3}\sum_{j=1}^n\frac{1}{j}\sum_{k=j}^n\frac{1}{n+k}\equiv\frac{H_n^2}{3}-\frac{1}{3}\sum_{k=1}^n\frac{H_{n+k-1}}{k}\pmod{p},
\end{align*}
and
\begin{align*}
\sum_{k=1}^n\frac{H_{n+k-1}}{k}=\sum_{k=1}^n\frac{H_{n+k}}{k}-\sum_{k=1}^n\frac{1}{k(n+k)}
\equiv \sum_{k=1}^n\frac{H_{n+k}}{k}+3H_n\pmod{p},
\end{align*}
and so
\begin{align}
\sum_{k=1}^n\frac{H_k}{3k-1}\equiv \frac{H_n^2}{3}-H_n-\frac{1}{3}\sum_{k=1}^n\frac{H_{n+k}}{k}\pmod{p}.\label{c-5}
\end{align}
Note that
\begin{align}
T_k=\frac{1}{3}\sum_{j=1}^k\frac{1}{j-2/3}\equiv \frac{1}{3}\sum_{j=1}^k\frac{1}{j+2n}
=\frac{1}{3}\left(H_{2n+k}-H_{2n}\right)\pmod{p}.\label{c-6}
\end{align}
By \eqref{c-6} and \eqref{b-1}, we have
\begin{align}
\sum_{k=1}^n\frac{T_k}{k}
\equiv \frac{1}{3}\sum_{k=1}^n\frac{H_{2n+k}-H_{2n}}{k}\equiv\frac{1}{3}\sum_{k=1}^n\frac{H_{2n+k}}{k}-\frac{H_n^2}{3}
\pmod{p}.\label{c-7}
\end{align}
Furthermore, by \eqref{c-6} we have
\begin{align*}
&\sum_{k=1}^n\frac{T_k}{3k-1}\equiv \frac{1}{3}\sum_{k=1}^n\frac{T_k}{n+k}\equiv \frac{1}{9}
\sum_{k=1}^n\frac{H_{2n+k}-H_{2n}}{n+k}
=\frac{1}{9}\sum_{k=n+1}^{2n}\frac{H_{n+k}}{k}-\frac{H_{2n}}{9}\sum_{k=1}^n\frac{1}{n+k}\pmod{p}.
\end{align*}
By \eqref{b-1}, we have
\begin{align*}
\sum_{k=1}^n\frac{1}{n+k}=H_{2n}-H_n\equiv 0\pmod{p}.
\end{align*}
It follows that
\begin{align}
\sum_{k=1}^n\frac{T_k}{3k-1}\equiv \frac{1}{9}\sum_{k=n+1}^{2n}\frac{H_{n+k}}{k}\pmod{p}.
\label{c-8}
\end{align}
Substituting \eqref{c-5}, \eqref{c-7} and \eqref{c-8} into the right-hand side of
\eqref{c-4} gives
\begin{align}
&\sum_{k=0}^{(p-1)/3}(6k+1)\frac{\left(\frac{1}{3}\right)_k^4 (1)_{2k}}
{(1)_k^4\left(\frac{2}{3}\right)_{2k}}\sum_{j=1}^k\left(\frac{1}{(3j)^2}-\frac{1}{(3j-2)^2}\right)\notag\\
&\equiv\frac{2p}{9}\left(\sum_{k=1}^n\frac{H_{2n+k}}{k}-\sum_{k=1}^{2n}\frac{H_{n+k}}{k}-\frac{9}{2}S_n^{(2)}\right)\pmod{p^2}.\label{c-9}
\end{align}

Using the same method as in Lemmas \ref{l-1} and \ref{l-2}, we can discover and prove the following two identities:
\begin{align}
&\sum_{k=1}^{2n}\frac{(-1)^k}{k}{3n\choose n+k}={3n\choose n}\left(H_n-H_{3n}\right),\label{c-10}\\
&\sum_{k=1}^{n}\frac{(-1)^k}{k}{3n\choose 2n+k}={3n\choose n}\left(H_{2n}-H_{3n}\right).\label{c-11}
\end{align}
Since for $0\le k\le p-1$,
\begin{align*}
{p-1\choose k}=\frac{(p-1)(p-2)\cdots(p-k)}{k!}\equiv (-1)^k-p(-1)^{k}H_k\pmod{p^2},
\end{align*}
and $n$ is an even integer, we have
\begin{align}
&{3n\choose n+k}\equiv (-1)^{k}-p(-1)^{k}H_{n+k}\pmod{p^2},\label{c-12}\\
&{3n\choose 2n+k}\equiv (-1)^{k}-p(-1)^{k}H_{2n+k}\pmod{p^2}.\label{c-13}
\end{align}
It follows from \eqref{c-10}--\eqref{c-13} that
\begin{align*}
&p\sum_{k=1}^{2n}\frac{H_{n+k}}{k}\equiv H_{2n}-{3n\choose n}H_n \pmod{p^2},\\
&p\sum_{k=1}^{n}\frac{H_{2n+k}}{k}\equiv H_{n}-{3n\choose n}H_{2n} \pmod{p^2},
\end{align*}
where we have utilized the fact that $H_{3n}\equiv 0\pmod{p^2}$.
Therefore, by \eqref{b-1} we have
\begin{align*}
p\left(\sum_{k=1}^{2n}\frac{H_{n+k}}{k}-\sum_{k=1}^{n}\frac{H_{2n+k}}{k}\right)
&\equiv (H_{2n}-H_n)\left(1+{3n\choose n}\right)\\
&\equiv pH_n^{(2)}\left(1+{3n\choose n}\right)\\
&\equiv 2pH_n^{(2)}\pmod{p^2},
\end{align*}
where we have used the fact that ${3n\choose n}\equiv 1\pmod{p}$. It follows that
\begin{align}
\sum_{k=1}^{2n}\frac{H_{n+k}}{k}-\sum_{k=1}^{n}\frac{H_{2n+k}}{k}\label{c-14}
\equiv 2H_n^{(2)}\pmod{p}.
\end{align}

Substituting \eqref{c-14} into \eqref{c-9} and using \eqref{b-4} and \eqref{new-1}, we obtain
\begin{align}
&\sum_{k=0}^{(p-1)/3}(6k+1)\frac{\left(\frac{1}{3}\right)_k^4 (1)_{2k}}
{(1)_k^4\left(\frac{2}{3}\right)_{2k}}\sum_{j=1}^k\left(\frac{1}{(3j)^2}-\frac{1}{(3j-2)^2}\right)\notag\\
&\equiv\frac{2p}{9}\left(2H_n^{(2)}-\frac{9}{2}S_n^{(2)}\right)\notag\\
&\equiv -\frac{p}{9}B_{p-2}\left(\frac{1}{3}\right)\pmod{p^2}.\label{c-15}
\end{align}
Then the proof of \eqref{a-3} follows from \eqref{c-3} and \eqref{c-15}.

\section{Proof of Theorem \ref{t-2}}
Let $n=(p+1)/3$. It is clear that $n$ is an even integer. Similarly to the proof of Theorem \ref{t-1}, by using \eqref{b-8} and \eqref{b-10} we can establish the following two supercongruences:
\begin{align}
&\sum_{k=0}^{(p+1)/3}(6k-1)\frac{\left(-\frac{1}{3}\right)_k^4 (1)_{2k}}
{(1)_k^4\left(-\frac{2}{3}\right)_{2k}}\notag\\
&\equiv p-p^2\sum_{k=0}^{(p+1)/3}(6k-1)\frac{\left(-\frac{1}{3}\right)_k^4 (1)_{2k}}
{(1)_k^4\left(-\frac{2}{3}\right)_{2k}}
\sum_{j=1}^k\left(\frac{1}{(3j)^2}-\frac{1}{(3j-4)^2}\right)
\pmod{p^4},\label{d-1}
\end{align}
and
\begin{align}
&\sum_{k=0}^{(p+1)/3}(6k-1)\frac{\left(-\frac{1}{3}\right)_k^4 (1)_{2k}}
{(1)_k^4\left(-\frac{2}{3}\right)_{2k}}
\sum_{j=1}^k\left(\frac{1}{(3j)^2}-\frac{1}{(3j-4)^2}\right)\notag\\
&\equiv \frac{p}{3}\left(H_n-3T_n^{(2)}+2\sum_{k=1}^n\frac{H_k}{3k-2}+2\sum_{k=1}^{n-1}\frac{S_k}{k}-6\sum_{k=1}^{n-1}\frac{S_k}{3k-2}\right)
\pmod{p^2}.\label{d-2}
\end{align}

Similarly to \eqref{c-5}, \eqref{c-7} and \eqref{c-8}, we can show that
\begin{align}
&\sum_{k=1}^n\frac{H_k}{3k-2}\equiv -\frac{1}{3}\sum_{k=n}^{2n-1}\frac{H_{n+k-1}}{k}-\frac{H_n}{2}\pmod{p},\label{d-3}\\
&\sum_{k=1}^{n-1}\frac{S_k}{k}\equiv -\frac{1}{3}H_{n-1}^2+\frac{1}{3}\sum_{k=1}^{n-1}\frac{H_{2n+k-1}}{k}\pmod{p},\label{d-4}\\
&\sum_{k=1}^{n-1}\frac{S_k}{3k-2}\equiv \frac{1}{9}\sum_{k=1}^{n-1}\frac{H_{n+k-1}}{k}
-\frac{H_n H_{n-1}}{9}+\frac{H_n}{3}\pmod{p}.\label{d-5}
\end{align}
Substituting \eqref{d-3}--\eqref{d-5} into the right-hand side of \eqref{d-2} gives
\begin{align}
&\sum_{k=0}^{(p+1)/3}(6k-1)\frac{\left(-\frac{1}{3}\right)_k^4 (1)_{2k}}
{(1)_k^4\left(-\frac{2}{3}\right)_{2k}}
\sum_{j=1}^k\left(\frac{1}{(3j)^2}-\frac{1}{(3j-4)^2}\right)\notag\\
&\equiv \frac{2p}{9}\left(
\sum_{k=1}^{n-1}\frac{H_{2n+k-1}}{k}-\sum_{k=1}^{2n-1}\frac{H_{n+k-1}}{k}
+H_nH_{n-1}-3H_n-H_{n-1}^2-\frac{9}{2}T_n^{(2)}\right)\label{d-6}
\pmod{p^2}.
\end{align}

Similarly to \eqref{c-14}, by using the following two identities discovered by {\tt Sigma}:
\begin{align*}
&\sum_{k=1}^{2n-1}\frac{(-1)^k}{k}{3n-2\choose n+k-1}={3n-2\choose n-1}\left(H_{n-1}-H_{3n-2}\right),\\
&\sum_{k=1}^{n-1}\frac{(-1)^k}{k}{3n-2\choose 2n+k-1}={3n-2\choose n-1}\left(H_{2n-1}-H_{3n-2}\right),
\end{align*}
we obtain
\begin{align*}
p\left(\sum_{k=1}^{n-1}\frac{H_{2n+k-1}}{k}-\sum_{k=1}^{2n-1}\frac{H_{n+k-1}}{k}\right)
&\equiv (H_{2n-1}-H_{n-1})\left({3n-2\choose n-1}-1\right)\pmod{p^2}.
\end{align*}
By \eqref{b-1} and ${3n-2\choose n-1}\equiv -1\pmod{p}$, we get
\begin{align}
\sum_{k=1}^{n-1}\frac{H_{2n+k-1}}{k}-\sum_{k=1}^{2n-1}\frac{H_{n+k-1}}{k}
&\equiv -2H_{n-1}^{(2)}\pmod{p}.\label{d-7}
\end{align}

Combining \eqref{d-6} and \eqref{d-7}, we arrive at
\begin{align}
&\sum_{k=0}^{n}(6k-1)\frac{\left(-\frac{1}{3}\right)_k^4 (1)_{2k}}
{(1)_k^4\left(-\frac{2}{3}\right)_{2k}}
\sum_{j=1}^k\left(\frac{1}{(3j)^2}-\frac{1}{(3j-4)^2}\right)\notag\\
&\equiv\frac{p}{3}\left(-3T_n^{(2)}-\frac{4}{3}H_{n-1}^{(2)}+\frac{2}{3}H_nH_{n-1}-2H_n-\frac{2}{3}H_{n-1}^2
\right)
\pmod{p^2}.\label{d-8}
\end{align}
Applying \eqref{new-1}, \eqref{b-6} and $H_n\equiv H_{n-1}+3\pmod{p}$ to the right-hand side
of \eqref{d-8} gives
\begin{align}
&\sum_{k=0}^{n}(6k-1)\frac{\left(-\frac{1}{3}\right)_k^4 (1)_{2k}}
{(1)_k^4\left(-\frac{2}{3}\right)_{2k}}
\sum_{j=1}^k\left(\frac{1}{(3j)^2}-\frac{1}{(3j-4)^2}\right)\notag\\
&\equiv p\left(\frac{1}{9}B_{p-2}\left(\frac{1}{3}\right)-2\right)
\pmod{p^2}.\label{d-9}
\end{align}
Then the proof of \eqref{a-4} follows from \eqref{d-1} and \eqref{d-9}.

\vskip 5mm \noindent{\bf Acknowledgments.}
This work was supported by the National Natural Science Foundation of China (grant 11801417).

\end{document}